\documentclass[12pt,oneside,english]{amsart}
\usepackage[T1]{fontenc}
\usepackage[latin1]{inputenc}
\usepackage{amssymb}

\makeatletter

 \theoremstyle{plain}
 \theoremstyle{plain}    
 \newtheorem{lem}{Lemma} 
 \theoremstyle{plain}    
 \newtheorem{thm}{Theorem} 
 \theoremstyle{plain}    
 \newtheorem{cor}{Corollary} 

\usepackage{babel}
\makeatother
\begin{document}

\title{A Free Analogue of Shannon's Problem on Monotonicity of Entropy}

\author{Dimitri Shlyakhtenko.}

\begin{abstract}
We prove a free probability analog of a result of \cite{shannon-monotonicity}.
In particualar we prove that if $X_{1},X_{2},\ldots$ are freely independent
identically distributed random variables, then the function\[
n\mapsto\chi\left(\frac{X_{1}+\cdots+X_{n}}{\sqrt{n}}\right)\]
is monotone increasing for all $n$.
\end{abstract}
\maketitle

\section{Introduction.}

Let $X$ be a random variable with law $\mu$ and let\[
\chi(X)=\int\int\log|s-t|d\mu(s)d\mu(t)+\frac{3}{4}+\frac{1}{2}\log2\pi.\]
This quantity, called free entropy, was discovered by Voiculescu in
\cite{dvv:entropy1} and plays the role of entropy in his free probability
theory (see e.g. \cite{dvv:entropysurvey} for a survey). Free entropy
has nice behavior with respect to freely independent random variables.
Amazingly, its behavior is in many instances parallel to the behavior
of classical entropy, if one replaces the classical notion of independence
by that of free independence. For example, for variables with variance
$1$, the free entropy is maximized by the semicircular law $\frac{1}{2\pi}\sqrt{2-t^{2}}\chi_{[-\sqrt{2},\sqrt{2}]}(t)dt$
(which plays the role of the Gaussian law in the free central limit
theorem) \cite{dvv:entropy4}. Similarly, one has the free analogue
of the entropy power inequality \cite{dvv-szarek:entropyPowerInequality}.

In analogy to the classical case, Voiculescu developed an infinitesimal
theory for $\chi$ with respect to free Brownian motion \cite{dvv:entropy5}.
If $S$ is a random variable with the semicircular law, then\[
\Phi(X)=\frac{d}{dt}\Big|_{t=0}\chi(X+\sqrt{t}S)\]
is called (by analogy with the classical setting) the free Fisher
information of $X$. If $X$ has law $\mu$ with density $d\mu(x)=f(x)dx$,
then\[
\Phi(X)=\frac{2}{3}\int(f(x))^{3}dx.\]
Free entropy can be recovered from the free Fisher information using
the formula\begin{equation}
\chi(X)=\frac{1}{2}\int\left(\frac{1}{1+t}-\Phi(X+\sqrt{t}S)\right)dt+\frac{1}{2}\log2\pi e\label{eq:chiInTermsOfPhi}\end{equation}
(a similar formula, which involves $\Phi(\sqrt{e^{2t}}X+\sqrt{1-e^{-2t}}S)$
can be derived by a change of variables). 

The free Fisher information also has many properties analogous to
those of classical Fisher information. These include the free analogs
of the Cramer-Rao inequality and also of the Stam inequality.

In \cite{shannon-monotonicity}, the authors have solved an old problem
going back to Shannon on monotonicity of entropy (see \cite{shannon-weaver:book,stam:inequalities,lieb:proofEntropyConjecture}).
In particular, they have proved that if $X_{1},X_{2},\ldots$ is a
sequence of iid random variables, then the classical entropy $H$
of their central limit sums is increasing:\[
H\left(\frac{X_{1}+\cdots+X_{n}}{\sqrt{n}}\right)\leq H\left(\frac{X_{1}+\cdots+X_{n+1}}{\sqrt{n+1}}\right).\]
Their proof relied on a new variational characterization of Fisher
information (Theorem 4 in \cite{shannon-monotonicity}).

The main purpose of this note is to derive the free analog of this
statement. Namely, if $X_{1},X_{2},\ldots$ are identically distributed
and freely independent, then the inequality\[
\chi\left(\frac{X_{1}+\cdots+X_{n}}{\sqrt{n}}\right)\leq\chi\left(\frac{X_{1}+\cdots+X_{n+1}}{\sqrt{n+1}}\right)\]
holds for all values of $n$.

The proof in the free case is in many instances a repetition of the
argument of \cite{shannon-monotonicity}, taken with minimal modifications;
we provide the full arguments for completeness. Indeed, once we obtain
a free analog of equation (4) of \cite{shannon-monotonicity} in Lemma
\ref{lemma:ineqForFisher} the rest of the argument is effectively
the same as in \cite{shannon-monotonicity}, with $H$ replaced by
$\chi$. However, we do not know of a free analog of the variational
characterization (Theorem 4 in \cite{shannon-monotonicity}), and
so we take a slightly different route to derive the analog of (4)
in \cite{shannon-monotonicity}, going back to the definition of the
free Fisher information using the {}``conjugate variables'' (which
are the free analogs of the classical score function), and is in its
spirit similar to the proof of the free Stam inequality \cite{dvv:entropy5}. 

It is curious to note that our proof can also be used, with appropriate
modifications, also in the classical case, thus giving a somewhat
shorter argument for the results of \cite{shannon-monotonicity},
avoiding the use of (Theorem 4 in \cite{shannon-monotonicity}). We
outline this approach in Section \ref{lemma:JisCondExpClassical}.
It would be interesting if this approach could be used to give a speedier
proof of the variational principle in \cite{shannon-monotonicity}.

As in \cite{shannon-monotonicity}, we prove results which are stronger
than monotonicity of entropy. These include generalizations to many
summands of the free Stam inequality and the free entropy power inequality.
These results are stated and proved in Section \ref{sec:Monotonicity-of-Free}.

\section{Monotonicity of Free Entropy.\label{sec:Monotonicity-of-Free}}

Let $(X,Y)\in(M,\tau)$ be two self-adjoint random variables. Following
\cite{dvv:entropy5}, we denote by $\partial_{X:Y}$ the derivation
from the algebra generated by $X$ and $Y$ into $L^{2}(M)\otimes L^{2}(M)$
determined by $\partial_{X:Y}(X)=1\otimes1$, $\partial_{X:Y}(Y)=0$
and the Leibniz rule. Then the conjugate variable $J(X:Y)$ (which
is the free analog of the classical score function) is the (unique,
if it exists) vector in $L^{2}(M,\tau)$ satisfying\[
\langle J(X:Y),P\rangle=\langle1\otimes1,\partial_{X:Y}(P)\rangle=\tau\otimes\tau(\partial_{X:Y}(P^{*})),\]
for any non-commuting polynomial $P$ in $X$ and $Y$. One write
$J(X)$ for $J(X:0)$. By definition, the free Fisher information
$\Phi(X)$ is given by\[
\Phi(X)=\Vert J(X)\Vert_{L^{2}(M)}^{2}.\]
We caution the reader familiar with \cite{shannon-monotonicity} that
in that paper the symbol $J$ denotes the classical analog of $\Phi$
and not of the conjugate variable $J$.

\begin{lem}
\label{lemma:JisCondExp}Let $X_{1},\ldots,X_{n}\in(M,\tau)$ be self-adjoint
random variables. Let $a_{1},\ldots,a_{n+1}\in\mathbb{R}$. Assume
that for some $j$, $J(\sum_{i\neq j}a_{i}X_{i}:X_{j})$ exists. Then\[
J(\sum_{i=1}^{n+1}a_{i}X_{i})=E_{W^{*}(\sum_{i=1}^{n+1}a_{i}X_{i})}(J(\sum_{i\neq j}a_{i}X_{i}:X_{j})).\]

\end{lem}
\begin{proof}
Let $f(t)=t^{m}$ be a monomial in one variable. Let $Y=\sum_{i=1}^{n+1}a_{i}X_{i}$,
$Y_{j}=\sum_{i\neq j}a_{i}X_{i}$, $Z_{j}=a_{j}X_{j}$, $N=W^{*}(Y)$.
Then\begin{eqnarray*}
\langle f(\sum_{i=1}^{n+1}a_{i}X_{i}),E_{N}(J(\sum_{i\neq j}a_{i}X_{i}:X_{j}))\rangle & = & \langle f(\sum_{i=1}^{n+1}a_{i}X_{i}),J(\sum_{i\neq j}a_{i}X_{i}:X_{j})\rangle\\
 & = & \tau\otimes\tau(\partial_{\sum_{i\neq j}a_{i}X_{i}:X_{j}}f(\sum_{i\neq j}a_{i}X_{i}+a_{j}X_{j}))\\
 & = & \tau\otimes\tau(\partial_{Y_{j}:Z_{j}}(Y_{j}+Z_{j})^{m})\\
 & = & \sum_{k=1}^{m}\tau\otimes\tau((Y_{j}+Z_{j})^{k-1}\otimes(Y_{j}+Z_{j})^{m-k})\\
 & = & \tau\otimes\tau(\partial_{Y}f(Y)).\end{eqnarray*}
Since $E_{N}(J(\sum_{i\neq j}a_{i}X_{i}:X_{j}))\in L^{2}(W^{*}(\sum_{i=1}^{n}a_{i}X_{i}$)),
the claimed equality follows.
\end{proof}
We get the following corollary:

\begin{lem}
\label{lemma:ineqForFisher}Let $X_{1},\ldots,X_{n+1}$ be freely
independent. Let $(a_{1},\ldots,a_{n+1})\in\mathbb{R}^{n+1}$ be an
$(n+1)$-tuple satisfying $\sum a_{j}^{2}=1$. Let $b_{1},\ldots,b_{n+1}\in\mathbb{R}$
be so that $\sum_{j=1}^{n+1}b_{j}\sqrt{1-a_{j}^{2}}=1$. Then\begin{equation}
\Phi(\sum_{j=1}^{n+1}a_{i}X_{i})\leq n\sum_{j=1}^{n+1}b_{j}^{2}\Phi\left(\frac{1}{\sqrt{1-a_{j}^{2}}}\sum_{i\neq j}a_{i}X_{i}\right).\label{eq:FisherInequality}\end{equation}

\end{lem}
\begin{proof}
Clearly, if $\Phi(\sum_{i\neq j}a_{i}X_{i})=+\infty$ for some $j$,
there is nothing to prove. Hence we assume that $\Phi(\sum_{i\neq j}a_{i}X_{i})$
is finite for all $j$.

Let $N=W^{*}(\sum_{i=1}^{n+1}a_{i}X_{i})$. According to Lemma \ref{lemma:JisCondExp},\[
J(\sum_{i=1}^{n+1}a_{i}X_{i})=E_{N}(J(\sum_{i\neq j}a_{i}X_{i}:X_{j})).\]
Hence\begin{eqnarray*}
J(\sum_{i=1}^{n+1}a_{i}X_{i}) & = & \sum_{j=1}^{n+1}b_{j}\sqrt{1-a_{j}^{2}}J(\sum_{i=1}^{n+1}a_{i}X_{i})\\
 & = & \sum_{j=1}^{n+1}b_{j}\sqrt{1-a_{j}^{2}}E_{N}J(\sum_{i\neq j}a_{i}X_{i}:X_{j})).\\
 & = & E_{N}(\sum_{j=1}^{n+1}b_{j}J(\frac{1}{\sqrt{1-a_{j}^{2}}}\sum_{i\neq j}a_{i}X_{i}:X_{j})).\end{eqnarray*}
Since $E_{N}$is a contraction on $L^{2}$ we deduce that\[
\Phi(\sum_{i=1}^{n+1}a_{i}X_{i})\leq\Vert\sum_{j=1}^{n+1}b_{j}J(\frac{1}{\sqrt{1-a_{j}^{2}}}\sum_{i\neq j}a_{i}X_{i}:X_{j})\Vert_{2}^{2}.\]
By freeness, $J(\frac{1}{\sqrt{1-a_{j}^{2}}}\sum_{i\neq j}a_{i}X_{i}:X_{j})=J(\frac{1}{\sqrt{1-a_{j}^{2}}}\sum_{i\neq j}a_{i}X_{i})$
and hence if we let\[
\xi_{j}=b_{j}J(\frac{1}{\sqrt{1-a_{j}^{2}}}\sum_{i\neq j}a_{i}X_{i}),\]
then\[
\Phi(\sum_{i=1}^{n+1}a_{i}X_{i})\leq\Vert\sum_{j=1}^{n+1}\xi_{j}\Vert_{2}^{2}.\]
Now let $E_{j}:M=W^{*}(X_{1},\ldots,X_{n})\to M_{j}=W^{*}(X_{1},\ldots,\hat{X}_{j},\ldots,X_{n})$
be the conditional expectation. Then $E_{j}:L^{2}(M)\to L^{2}(M)$
are projections and moreover $E_{j}$ form a commuting family. Indeed,
because of the freeness assumptions, we may write\begin{eqnarray*}
M & = & W^{*}(X_{1})*\cdots*W^{*}(X_{j})*\cdots*W^{*}(X_{n})\\
E_{j} & = & \textrm{id}*\cdots*\tau|_{W^{*}(X_{j})}*\cdots*\textrm{id}.\end{eqnarray*}
Hence if $i<j$,\[
E_{i}\circ E_{j}=\textrm{id}*\cdots*\tau|_{W^{*}(X_{i})}*\textrm{id}*\cdots*\textrm{id}*\tau|_{W^{*}(X_{j})}*\textrm{id}*\cdots*\textrm{id}=E_{j}\circ E_{i}.\]
In particular, note that $E_{1}\circ\cdots\circ E_{n}=\tau$. Since
$\tau(\xi_{j})=0$ (because $\xi_{j}$ is up to a multiple a conjugate
variable), we deduce that conditions of Lemma 5 on p. 6 of \cite{shannon-monotonicity}
are satisfied. We apply this lemma to conclude that\[
\Phi(\sum_{i=1}^{n+1}a_{i}X_{i})\leq\Vert\sum_{j=1}^{n+1}\xi_{j}\Vert_{2}^{2}\leq n(\sum_{j=1}^{n+1}\Vert\xi_{j}\Vert_{2}^{2}).\]
Recalling that\[
\xi_{j}=b_{j}J(\frac{1}{\sqrt{1-a_{j}^{2}}}\sum_{i\neq j}a_{i}X_{i})\]
we conclude that\[
\Phi(\sum_{i=1}^{n+1}a_{i}X_{i})\leq n(\sum_{j=1}^{n+1}\Vert\xi_{j}\Vert_{2}^{2})=n\sum_{j=1}^{n+1}b_{j}^{2}\Phi\left(\frac{1}{\sqrt{1-a_{j}^{2}}}\sum_{i\neq j}a_{i}X_{i}\right),\]
as claimed.
\end{proof}
We can now deduce a many-variable version of the free Stam inequality
\cite{dvv:entropy5}: 

\begin{thm}
\label{theorem:freeStam}Let $X_{1},\ldots,X_{n+1}$ be free random
variables. Then\[
\frac{n}{\Phi(\sum_{i=1}^{n+1}X_{i})}\geq\sum_{j=1}^{n+1}\frac{1}{\Phi(\sum_{i\neq j}X_{i})}\]

\end{thm}
\begin{proof}
Let $a_{j}=(n+1)^{-1/2}$. Then $1/\sqrt{1-a_{j}^{2}}=\sqrt{(n+1)/n}$
and hence (\ref{eq:FisherInequality}) implies that for any $b_{j}$
with $\sum b_{j}=\sqrt{(n+1)/n},$ we have\[
\Phi(\sum_{i=1}^{n+1}X_{i})\leq n\sum_{j=1}^{n+1}\frac{n}{n+1}\  b_{j}^{2}\Phi(\sum_{i\neq j}X_{i}).\]
Hence if we are given $\lambda_{j}$ with $\sum\lambda_{j}=1$, we
could take $b_{j}=\lambda_{j}\sqrt{(n+1)/n}$ and hence deduce that\begin{equation}
\Phi(\sum_{i=1}^{n+1}X_{i})\leq n\sum_{j=1}^{n+1}\lambda_{j}^{2}\Phi(\sum_{i\neq j}X_{i}).\label{eq:lambdasFisher}\end{equation}
Now let $C=\sum_{j=1}^{n+1}\Phi(\sum_{i\neq j}X_{i})^{-1}$ and let\[
\lambda_{j}=\frac{\Phi(\sum_{i\neq j}X_{i})^{-1}}{C}.\]
Then (\ref{eq:lambdasFisher}) becomes\[
\Phi(\sum_{i=1}^{n+1}X_{i})\leq n\sum_{j=1}^{n+1}\frac{1}{C^{2}}\Phi(\sum_{i\neq j}X_{i})^{-1}=nC^{-1}.\]
Recalling the definition of $C$ gives\[
\Phi(\sum_{i=1}^{n+1}X_{i})\leq n\sum_{j=1}^{n+1}\Phi(\sum_{i\neq j}X_{i})^{-1},\]
which clearly implies the statement of the Theorem.
\end{proof}
\begin{thm}
\label{theorem:betterMonotonicity}Let $X_{1},\ldots,X_{n+1}$ be
free random variables and let $(a_{1},\ldots,a_{n+1})\in S^{n}$ be
a unit vector. Then\begin{equation}
\chi(\sum_{i=1}^{n+1}a_{i}X_{i})\geq\sum_{j=1}^{n+1}\frac{1-a_{j}^{2}}{n}\chi\left(\frac{1}{\sqrt{1-a_{j}^{2}}}\sum_{i\neq j}a_{i}X_{i}\right).\label{eq:chiIneq}\end{equation}

\end{thm}
\begin{proof}
Set $b_{j}=\frac{1}{n}\sqrt{1-a_{j}^{2}}$ in Lemma \ref{lemma:ineqForFisher}
and apply the Lemma to $X_{j}^{(t)}=X_{j}+\sqrt{t}S_{j}$. We then
get\[
\Phi(\sum_{i=1}^{n+1}a_{i}X_{i}^{(t)})\leq\sum_{j}\frac{1-a_{j}^{2}}{n}\Phi\left(\frac{1}{\sqrt{1-a_{j}^{2}}}\sum_{i\neq j}a_{i}X_{i}^{(t)}\right).\]
Notice that since $\sum a_{i}^{2}=1$, the law of $\sum_{i=1}^{n+1}a_{i}X_{i}^{(t)}$
is the same as the law of $\sum_{i=1}^{n+1}a_{i}X_{i}+\sqrt{t}S'$,
where $S'=\sum_{i=1}^{n+1}a_{j}S_{j}$ is a semicircular variable,
free from $\sum_{i=1}^{n+1}a_{i}X_{i}$. Similarly, since $\sum_{i\neq j}a_{i}^{2}=(1-a_{j}^{2})$,
we get that the law of $\frac{1}{\sqrt{1-a_{j}^{2}}}\sum_{i\neq j}a_{i}X_{i}^{(t)}$
is the same as the law of $(\frac{1}{\sqrt{1-a_{j}^{2}}}\sum_{i\neq j}a_{i}X_{i})+\sqrt{t}S^{(j)}$,
where $S^{(j)}=\frac{1}{\sqrt{1-a_{j}^{2}}}\sum_{i\neq j}a_{j}S_{j}$
is a semicircular free from $\sum_{i\neq j}a_{i}X_{i}$. Hence using
the formula (\ref{eq:chiInTermsOfPhi}) and integrating with respect
to $t$, we deduce the desired inequality (\ref{eq:chiIneq}).
\end{proof}
\begin{cor}
\label{corrollary:monotonicity}Let $X_{1},\ldots,X_{n+1}$ be identically
distributed free random variables. Then\[
\chi\left(\frac{X_{1}+\cdots+X_{n}}{\sqrt{n}}\right)\leq\chi\left(\frac{X_{1}+\cdots+X_{n+1}}{\sqrt{n+1}}\right).\]

\end{cor}
\begin{proof}
Using Theorem \ref{theorem:betterMonotonicity} with $a_{j}=1/\sqrt{n+1}$
gives us\begin{eqnarray*}
\chi(\frac{1}{\sqrt{n+1}}\sum_{i=1}^{n+1}X_{i}) & \geq & \sum_{j=1}^{n+1}\frac{1-\frac{1}{n+1}}{n}\chi\left(\frac{1}{\sqrt{1-\frac{1}{n+1}}}\sum_{i\neq j}\frac{1}{\sqrt{n+1}}X_{i}\right)\\
 & = & \sum_{j=1}^{n+1}\frac{1}{n+1}\chi(\frac{1}{\sqrt{n}}\sum_{i\neq j}X_{i})\\
 & = & \chi(\frac{1}{\sqrt{n}}\sum_{i\neq j}X_{i}),\end{eqnarray*}
since $X_{1},\ldots,X_{n}$ have the same law.
\end{proof}
One can also get a free analogue of the entropy power inequality for
many summands:

\begin{thm}
\label{theorem:entropyPower}Let $X_{1},\ldots,X_{n+1}$ be free random
variables. Then\[
\exp\left[2\chi(\sum_{i=1}^{n+1}X_{i})\right]\geq\frac{1}{n}\sum_{j=1}^{n+1}\exp\left[2\chi(\sum_{i\neq j}X_{i})\right].\]

\end{thm}
\begin{proof}
Let $E=\exp(2\chi(\sum_{i=1}^{n+1}X_{i}))$ and $E_{j}=\exp(2\chi(\sum_{i\neq j}X_{i}))$.
Note that if $X$ and $Y$ are free, then using the triangular change
of variables formula \cite{dvv:entropy2} and monotonicity of free
entropy, one gets:\[
\chi(X+Y)+\chi(Y)\geq\chi(X+Y,Y)=\chi(X,Y)=\chi(X)+\chi(Y),\]
so that $\chi(X+Y)\geq\chi(X)$. Hence $E\geq E_{j}$ for all $j$.
Hence if for some $j$, $E_{j}\geq\frac{1}{n}\sum_{i=1}^{n+1}E_{i}$,
the desired inequality $E\geq\frac{1}{n}\sum E_{j}$ would follow.
Thus we may assume that $E_{j}<\frac{1}{n}\sum_{i=1}^{n+1}E_{i}$.
Let then $\lambda_{j}=E_{j}/\sum_{i=1}^{n+1}E_{i}$; thus $\lambda_{j}<1/n$
for all $j$. Let $a_{j}=\sqrt{1-n\lambda_{j}}$; applying (\ref{eq:chiIneq})
gives the inequality\begin{equation}
\chi(\sum_{i=1}^{n+1}X_{i})\geq\sum_{j=1}^{n+1}\lambda_{j}\chi(\frac{1}{\sqrt{n\lambda_{j}}}\sum_{i\neq j}X_{i})\label{eq:chiInsidePower}\end{equation}
Using the fact that\begin{eqnarray*}
\chi(\frac{1}{\sqrt{n\lambda_{j}}}\sum_{i\neq j}X_{i}) & = & \chi(\sum_{i\neq j}X_{i})-\frac{1}{2}\log n\lambda_{j}\\
 & = & \chi(\sum_{i\neq j}X_{i})-\chi(\sum_{i\neq j}X_{i})+\frac{1}{2}\log\frac{1}{n}\sum_{i=1}^{n+1}E_{i}\\
 & = & \frac{1}{2}\log(\frac{1}{n}\sum_{i=1}^{n+1}E_{i}),\end{eqnarray*}
(\ref{eq:chiInsidePower}) becomes\[
\chi(\sum_{i=1}^{n+1}X_{i})\geq\sum_{j=1}^{n+1}\lambda_{j}\frac{1}{2}\log(\frac{1}{n}\sum_{i=1}^{n+1}E_{i})=\frac{1}{2}\log(\frac{1}{n}\sum_{i=1}^{n+1}E_{i}).\]
Multiplying by two and exponentiating gives\[
E\geq\frac{1}{n}\sum_{i=1}^{n+1}E_{i},\]
which is the desired inequality.\[
\]

\end{proof}

\section{\label{sec:Some-Remarks-on}Some Remarks on the Classical Case.}

The argument in the free case followed almost word for word the argument
of \cite{shannon-monotonicity} in the classical case with one exception:
we used Lemma \ref{lemma:JisCondExp} to prove Lemma \ref{lemma:ineqForFisher}
(which is the free analog of equation (4) in \cite{shannon-monotonicity}).
Once this analog of (4) in \cite{shannon-monotonicity} was established,
the argument in the free case became entirely parallel to the classical
case.

We now point out that an analog of Lemma \ref{lemma:JisCondExp} also
holds in the classical case, and hence equation (4) in \cite{shannon-monotonicity}
can be derived this way (at least for variables that have finite moments
of all orders), avoiding the use of the variational characterization
of classical Fisher information (Theorem 4 in \cite{shannon-monotonicity}).
This approach seems to be slightly shorter and may shed more light
at the variational characterization mentioned above.

In the remainder of the section we'll recall some facts about the
classical Fisher information and indicate how the analog of Lemma
\ref{lemma:JisCondExp} can be proved. In order to avoid confusion
with the notation in the free case, we shall use $F$ to denote Fisher
information.

Let $X,Y_{1},\ldots,Y_{n}$ be classical real-valued random variables,
and let $\omega$ be the probability measure on $\mathbb{R}^{n+1}$
describing their joint law. Recall that the score function\[
j(X:Y_{1},\ldots,Y_{n})\in L^{2}(\omega)\]
is the (unique, if it exists) element satisfying\[
\langle j(X:Y_{1},\ldots,Y_{n}),f(x,y_{1},\ldots,y_{n})\rangle=\langle1,\frac{\partial f}{\partial x}\rangle.\]
In other words, $j(X:Y_{1},\ldots,Y_{n})=\left(\frac{\partial}{\partial x}\right)^{*}1$
if one considers $\partial/\partial x$ to be a densely defined operator
on $L^{2}(\omega)$ with domain of definition consisting of polynomials
in the coordinates $x,y_{1},\ldots,y_{n}$. Algebraically, $\partial/\partial x$
is determined by the fact that it is a derivation and also by $(\partial/\partial x)x=1$,
$(\partial/\partial x)y_{k}=0$ for all $k=1,\ldots,n$. When this
algebraic definition is emphasized, we'll denote $\partial/\partial x$
by $d_{x:y_{1},\ldots,y_{n}}$.$ $

In the case that $ $$n=0$, one can easily check that if $d\omega(x)=f(x)dx$,
then\[
j(X)=\frac{f'}{f}\]
and hence the Fisher information is given by\[
F(X)=\int\frac{(f'(x))^{2}}{f(x)}dx=\int\frac{(f'(x))^{2}}{f(x)^{2}}f(x)dx=\Vert j(X)\Vert_{L^{2}(\omega)}^{2}.\]

\begin{lem}
\label{lemma:JisCondExpClassical}Let $a_{1},\ldots,a_{n+1}\in\mathbb{R}$.
Assume that for some $k$, $j(\sum_{i\neq k}a_{i}X_{i}:X_{k})$ exists.
Then\[
j(\sum_{i=1}^{n+1}a_{i}X_{i})=E_{W^{*}(\sum_{i=1}^{n+1}a_{i}X_{i})}(j(\sum_{i\neq k}a_{i}X_{i}:X_{k})).\]

\end{lem}
\begin{proof}
Let $f(t)=t^{m}$ be a monomial in one variable. Let $Y=\sum_{i=1}^{n+1}a_{i}X_{i}$,
$Y_{k}=\sum_{i\neq k}a_{i}X_{i}$, $Z_{k}=a_{k}X_{k}$, $N=W^{*}(Y)$.
Then\begin{eqnarray*}
\langle f(\sum_{i=1}^{n+1}a_{i}X_{i}),E_{N}(j(\sum_{i\neq k}a_{i}X_{i}:X_{k}))\rangle & = & \langle f(\sum_{i=1}^{n+1}a_{i}X_{i}),j(\sum_{i\neq k}a_{i}X_{i}:X_{k})\rangle\\
 & = & \langle f(Y_{k}+Z_{k}),j(Y_{k}:Z_{k})\rangle\\
 & = & \langle d_{Y_{k}:Z_{k}}f(Y_{k}+Z_{k}),1\rangle\\
 & = & \sum_{r=1}^{m}\langle(Y_{k}+Z)^{r-1}\ (d_{Y_{k}:Z_{k}}(Y_{k}+Z_{k}))\ (Y_{k}+Z)^{m-r},1\rangle\\
 & = & \langle f'(Y),1\rangle.\end{eqnarray*}
Since $E_{N}(j(\sum_{i\neq k}a_{i}X_{i}:X_{k}))\in L^{2}(W^{*}(\sum_{i=1}^{n}a_{i}X_{i}$)),
the claimed equality follows.
\end{proof}
Finally, it is routine to check that $j(X:Y)=j(X)$ in the case that
$X$ and $Y$ are classically independent. Thus the argument in Lemma
\ref{lemma:ineqForFisher} applies verbatim to yield the classical
analog of (\ref{eq:FisherInequality}) and thus of equation (4) in
\cite{shannon-monotonicity}.

\bibliographystyle{amsalpha}

\begin{thebibliography}{ABBN04}

\bibitem[ABBN04]{shannon-monotonicity}
S.~Artstein, K.~Bally, F.~Barthez, and A.~Naor, \emph{Solution of {Shannon's}
  problem on monotonicity of entropy}, Journal Amer. Math. Soc. \textbf{17}
  (2004), 975--982.

\bibitem[Lie78]{lieb:proofEntropyConjecture}
E.H. Lieb, \emph{Proof of an entropy conjecture of {Wehrl}}, Comm. Math. Phys.
  \textbf{62} (1978), 35--41.

\bibitem[Sta59]{stam:inequalities}
A.~Stam, \emph{Some inequalities satisfied by the quantities of information of
  {Fisher} and {Shannon}}, Info. control \textbf{2} (1959), 101--112.

\bibitem[SV96]{dvv-szarek:entropyPowerInequality}
S.~Szarek and D.~Voiculescu, \emph{Volumes of restricted {Minkowski} sums and
  the free analogue of the entropy power inequality}, Comm. Math. Phys.
  \textbf{178} (1996), 563--570.

\bibitem[SW49]{shannon-weaver:book}
C.~Shannon and W.~Weaver, \emph{The mathematical theory of communication},
  Univ. Illinois Press, Urbana, IL, 1949.

\bibitem[Voi93]{dvv:entropy1}
D.-V. Voiculescu, \emph{The analogues of entropy and of {Fisher's} information
  measure in free probability theory {I}}, Commun. Math. Phys. \textbf{155}
  (1993), 71--92.

\bibitem[Voi94]{dvv:entropy2}
D.-V. Voiculescu, \emph{The analogues of entropy and of {Fisher's} information
  measure in free probability theory {II}}, Invent. Math. \textbf{118} (1994),
  411--440.

\bibitem[Voi97]{dvv:entropy4}
D.-V. Voiculescu, \emph{The analogues of entropy and of {Fisher}'s information
  measure in free probability theory, {IV}: Maximum entropy and freeness}, Free
  Probability (D.-V. Voiculescu, ed.), American Mathematical Society, 1997,
  pp.~293--302.

\bibitem[Voi98]{dvv:entropy5}
D.-V. Voiculescu, \emph{The analogues of entropy and of {Fisher}'s information
  measure in free probability, {V}}, Invent. Math. \textbf{132} (1998),
  189--227.

\bibitem[Voi02]{dvv:entropysurvey}
D.-V. Voiculescu, \emph{Free entropy}, Bull. London Math. Soc. \textbf{34}
  (2002), no.~3, 257--278.

\end{thebibliography}

\providecommand{\bysame}{\leavevmode\hbox to3em{\hrulefill}\thinspace}

\end{document}